\documentclass{article}
\usepackage{amsmath,amssymb}

\newcommand{\TeXmacs}{T\kern-.1667em\lower.5ex\hbox{E}\kern-.125emX\kern-.1em\lower.5ex\hbox{\textsc{m\kern-.05ema\kern-.125emc\kern-.05ems}}}
\newcommand{\email}[1]{{\textit{Email:} \texttt{#1}}}
\newcommand{\homepage}[1]{{\textit{Web:} \texttt{#1}}}
\newcommand{\tmop}[1]{\ensuremath{\operatorname{#1}}}
\newcommand{\tmtextit}[1]{{\itshape{#1}}}
\newcommand{\tmtextsl}[1]{{\slshape{#1}}}
\newcommand{\tmtextup}[1]{{\upshape{#1}}}
\newcommand{\withTeXmacstext}{This document has been produced using \TeXmacs (see \texttt{http://www.texmacs.org})}
\newtheorem{corollary}{Corollary}

\begin{document}

\title{You Could Simplify Calculus\thanks{{\withTeXmacstext}}}\author{Michael
Livshits\thanks{\email{michaelliv@gmail.com; };
\homepage{http://www.mathfoolery.org}}\\
Cambridge, MA, USA 02138 }\maketitle

\begin{abstract}
  I explain a direct approach to differentiation and integration. Instead of
  relying on the general notions of real numbers, limits and continuity, we
  treat functions as the primary objects of our theory, and view
  differentiation as division of $f (x) - f (a)$ by $x - a$ in a certain class
  of functions. When $f$ is a polynomial, the division can be carried out
  explicitly. To see why a polynomial with a positive derivative is increasing
  (the monotonicity theorem), we use the estimate $|f (x) - f (a) - f' (a) (x
  - a) | \leqslant K (x - a)^2 $. By making it into a definition we arrive at
  the notion of uniform Lipschitz differentiability (ULD), and see that the
  derivative of a ULD function is Lipschitz. Taking different moduli of
  continuity instead of the absolute value, we get different flavors of
  calculus, each rather elementary, but all together covering the total range
  of uniformly differentiable functions. Using the class of functions
  continuous at $a$, we recapture the classical notion of pointwise
  differentiability. It turns out that uniform Lipschitz differentiability is
  equivalent to divisibility of $f (x) - f (a)$ by $x - a$ in the class of
  Lipschitz functions of two variables, $x$ and $a.$ The same is true for any
  subadditive modulus of continuity. In this bottom-up, computational, one
  modulus of continuity at a time approach to calculus, the monotonicity
  theorem takes the central stage and provides the aspects of the subject that
  are important for practical applications. The weighty ontological issues of
  compactness and completeness can be treated lightly or postponed, since they
  are hardly used this streamlined approach that pretty much follows the
  Vladimir Arnold's ``principle of minimal generality, according to which
  every idea should first be understood in the simplest situation; only then
  can the method developed be applied to more complicated cases.'' I discuss a
  generalization to many variables briefly.
\end{abstract}

\begin{center}
  \section{Two Stories, One Fictional, One Real}
\end{center}

\subsection{Differentiating $x^4$ without using limits}

A teacher \tmtextit{}asks a student \tmtextsl{}to calculate the derivative of
$x^4$ at $x = a$. The student writes down the difference quotient $\frac{x^4 -
a^4}{x - a}$ , then, by factoring the numerator, rewrites it as $\frac{(x - a)
(x + a) (x^2 + a^2)}{x - a}$ , then cancels $x - a$ and gets $x + a$, then
substitutes $x = a$ and gets $4 a^3$, that is the right answer, of course. The
teacher does not like the solution, and the following conversation takes
place.

T: Your answer is correct, but why didn't you use the definition of the
derivative as a limit? We are studying calculus here, you know.

S: Do I really need to use limits? It looks like a waste of time, I can just
simplify and plug in $x = a$ instead, it looks like it works fine.

T: But do you understand why it works?

S: Hmmm, let me see. I guess it works because the limit of ($x + a) (x^2 +
a^2)$ as $x \rightarrow a$ is $4 a^3$, so, instead of calculating the limit we
can just plug $x = a$ into ($x + a) (x^2 + a^2)$.

T: How do you call such a function, that you can just plug in $x = a$ into it
instead of calculating the limit of this function at $a ?$

S: Continuous at $a$? Yeah, I remember.

T: Right! You know, people differentiated polynomials, roots and trig
functions in the 17th century, long before they started thinking of such
generalities as continuity and limits in the 19th century. Why don't you try
to differentiate your way some other simple algebraic expressions, such as
$\sqrt[3]{x}$ or $\frac{x^2}{3 + x^3}$ ?

S: O.K., I will, I think I understand it a little better now.

\subsection{Differentiating $\sqrt{x}$ without using limits}

It happened in the fall of 1997, when I taught two calculus recitation
sections at Suffolk University. The purpose of these sections was to answer
the questions the students had about their homework and the subject in
general. The text was Anton's Calculus, which I came to hate as the semester
progressed.

It was one of the classes, and some students asked me to explain how to
differentiate $\sqrt{x}$. So I wrote down the difference quotient
$\frac{\sqrt{x} - \sqrt{a}}{x - a} $on the chalkboard and said that we had to
calculate the limit of this expression as $x$ approaches $a$.

As soon as I uttered the word ``limit'' I saw many students slump in their
seats, their eyes glazing over, and I had the sinking feeling that they were
totally lost. I had to do something fast to help them, to pull them out of
their despair, but what?

I said, look, you don't really need limits to calculate this derivative, you
can do it algebraically. Let us rewrite this expression in such a way that it
would make sense for $x = a$. How can we do that? Let us rewrite the
denominator as $\sqrt{x}^2 - \sqrt{a}^2$ and factor it as $( \sqrt{x} -
\sqrt{a}) ( \sqrt{x} + \sqrt{a})$, so $\text{ $\frac{\sqrt{x} - \sqrt{a}}{x -
a} $} = \frac{\sqrt{x} - \sqrt{a}}{\text{$( \sqrt{x} - \sqrt{a}) ( \sqrt{x} +
\sqrt{a})$}} = \frac{1}{\sqrt{x} + \sqrt{a}}$ that makes sense for $x = a$,
giving us the answer ($\sqrt{x})' = 1 / (2 \sqrt{x})$, that's all there is to
it.

I saw the students brightening up a little bit, when they realized that the
problem could be solved with the tools familiar to them. And that's exactly
when it dawned on me that all calculus could be done like that,
differentiation being nothing but division in the class of continuous
functions. It surely looked like a promising idea.

\section{Calculus of Polynomials}

\subsection{Formal differentiation}

Let us start with the simplest and most popular example, differentiating
$x^2$. We form the difference quotient $\frac{x^2 - a^2}{x - a} $ and try to
make sense of it for $x = a$. The trouble is, of course, that when we just
plug in $x = a$, we get $0 / 0$, which is undefined, because $0 c = 0$ for any
number $c$. But luckily, the numerator factors as $(x - a) (x + a)$, so we can
cancel $x - a \tmop{and} \tmop{rewrite} \tmop{our} \tmop{expression} \tmop{as}
x + a$ that makes sense for $x = a$, giving us $(x^2)' = 2 x.$ To generalize
to $x^n$, we use the factorization $x^n - a^n = (x - a) (x^{n - 1} + x^{n - 2}
a + \ldots . + a^{n - 1})$ to get $(x^n)' = n x^{n - 1}$.

This trick will work for any polynomial $p (x)$, because $a$ is a root of the
polynomial $p (x) - p (a)$, and therefore $\tmop{it} \tmop{is}
\tmop{divisible} \tmop{by} x - a$, so we have $p (x) - p (a) = (x - a) q (x,
a),$ and we can rewrite $\frac{p (x) - p (a)}{x - a}$ as $q (x, a)$ which is a
polynomial in $x$ and $a$ and therefore makes sense for $x = a$, giving us $p'
(x) .$

Of course we don't have to divide polynomials every time we differentiate
them. The first two \tmtextsl{differentiation rules\tmtextup{}} tell us that
$(f + g)' = f' + g'$ and $(k f)' = k f' $ for any constant $k$, in other
words, differentiation is a linear operation, and therefore we can
differentiate polynomials ``term by term,'' i.e.
\[ (p_0 + p_1 x + \ldots + p_n x^n)' = p_1 + 2 p_2 x + \ldots + n p_n x^{n -
   1} . \]
The other two rules of differentiation, the \tmtextsl{product} (or
\tmtextsl{Leibniz}) rule, saying that $(f g)' = f' g + f g' $ and the
\tmtextsl{chain} rule by Newton, $(f (g (x)))' = f' (g (x)) g' (x)$ are \ a
matter of algebra of polynomials.

The trick developed here can be used to differentiate all rational functions,
and even algebraic functions that are defined implicitly by algebraic
equations, if we use implicit differentiation. \

\subsection{Double roots and the basic estimate}

Consider a polynomial $p (x) .$The question is: ``why the \tmtextsl{tangent
}to the graph $y = p (x)$ at the point $(a, p (a))$, which is the line defined
by the equation $y = p (a) + p' (a) (x - a)$ looks like a tangent, i.e.
``clings'' to this graph?'' \ Let us start with a simple example, $p (x) =
x^k$. Then $p' (a) = k a^{k - 1}$, and $x^k - a^k - k a^{k - 1} (x - a) = (x -
a) (x^{k - 1} + x^{k - 2} a + \ldots + a^{k - 1} - k a^{k - 1}) = (x - a)^2 r
(x, a)$, with $r$ a polynomial in \ $x$ and $a$, because the second factor
vanishes for $x = a$, so it is divisible by $x - a$. A similar factoring, $p
(x) - p (a) - p' (a) (x - a) = (x - a)^2 r (x, a)$, holds for any polynomial
$p$ since it is a sum of monomials. It shows that $x = a$ is a double root of
the equation $p (x) - p (a) - p' (a) (x - a) = 0$. This fact can be taken as
the definition of a tangent to a graph of a polynomial, and can be used to
define the derivative for polynomials. The vertical distance $d (x, a)$
between the graph and the tangent can be written as $(x - a)^2 |r (x, a) |$,
with $r$ a polynomial in \ $x$ and $a$. When $x$ and $a$ are contained in some
finite interval, \ $|r (x, a) |$ will be bounded from above by some constant
$K$, giving us an estimate $d (x, a) \leqslant K (x - a)^2 .$ This
\tmtextsl{basic estimate}, that also can be written as
\begin{equation}
  |p (x) - p (a) - p' (a) (x - a) | \leqslant K (x - a)^2
\end{equation}
holds for any polynomial $p,$ and explains why tangents clings the graphs. We
will use it in the next subsection to understand why a polynomial with a
positive derivative is increasing.

\subsection{Monotonicity principle}

The derivative is a mathematical metaphor for the instantaneous velocity, or
the instantaneous rate of change of a function relative to its argument. So we
would expect that a function with positive derivative would be increasing. Let
us see why it is true for polynomials. Assume that $p' (x) \geqslant 0$ for
any $x$ such that $A \leqslant x \leqslant B.$ We want to show that $p (A)
\leqslant p (B) .$ We can deal with a simpler case $p' (x) \geqslant C > 0$
first. Our basic estimate (1) tells us that $p (x) - p (a) \geqslant p' (a) (x
- a) - K (x - a)^2,$ so $p (a) \leqslant p (x)$ if $0 < x - a \leqslant C /
K$. Therefore, $p (A) \leqslant p (B),$ since we can get from $A$to $B$ by
steps shorter than $C / K$. To get to the original assumption, we can consider
$q (x) = p (x) + C x$ with $C > 0$ and conclude that $p (B) - p (A) \geqslant
C (A - B)$, therefore $p (A) \leqslant p (B)$ since $C$ is arbitrary.

By applying our monotonicity principle to $f + M x$ and $f - M x$, we can
demonstrate

\begin{corollary}
  The Rule of Bounded Change.
\end{corollary}

If $|p' | \leqslant M$, then $|p (x) - p (a) | \leqslant M|x - a|$.

When we look at definite integrals as increments of anti-derivatives, we can
see how monotonicity is related to positivity of the area.

\subsection{Formal integration}

It can be introduced before the basic estimate is treated and monotonicity
theorem is demonstrated, and it is very easy for polynomials. Besides, it
provides a strong evidence for the Newton-Leibniz theorem. The simplest
examples of course are the constants and the linear functions. A bit more work
is required to calculate the areas under the other power curves, and may give
the skeptics an opportunity to use such tools as algebra, the geometric series
and even combinatorics (to estimate the sum $1^k + 2^k + \ldots + n^k) .$
Newton-Leibniz is very intuitive and can be explained early on. The
integration rules are just the rules of differentiation, rewritten in terms of
integrals. This formal theory can be used right away to solve some interesting
problems in geometry and physics.

\section{Uniform Lipschitz Calculus}

How can we extend our calculus to functions more general than polynomials? As
it often happens in mathematics, we just look at some useful property or a
formula and make it into a definition (think about the Pythagorean Theorem).
The useful property here will be the basic estimate (1) from section 2.2, so
we call a function $f$ \tmtextsl{uniformly Lipschitz differentiable} (ULD) if
the estimate
\begin{equation}
  |f (x) - f (a) - f' (a) (x - a) | \leqslant K (x - a)^2
\end{equation}
holds for some constant $K \tmop{independent} \tmop{of}$ $x$ and $a$.

Now we can prove our monotonicity theorem from section 2.3 for ULD functions.

\subsection{The automatic Lipschitz estimate for the derivative}

We know that the derivatives of polynomials are also polynomials. What would
be the analogous fact for ULD functions? It turns out that theis derivatives
are \tmtextsl{Lipschitz}, i.e., they satisfy the estimate $|f' (x) - f' (a) |
\leqslant L|x - a|$ with $L$ independent of $x$ and $a$.

To see it, we notice that for $x \neq a$ $| \frac{f (x) - f (a)}{x - a} - f'
(a) | \leqslant K|x - a|$. By interchanging $x$ and $a$ we get $| \frac{f (a)
- f (x)}{a - x} - f' (x) | \leqslant K|a - x|$. but $\frac{f (a) - f (x)}{a -
x} = \frac{f (x) - f (a)}{x - a}$, and we see that $|f' (x) - f' (a) |
\leqslant 2 K|x - a|$, i.e., $f'$ is Lipschitz.

Of course all the polynomials are Lipschitz on any finite interval, because $x
- a$ is a factor in $p (x) - p (a)$, \ and the ULD functions are too, because
their derivatives are bounded on any finite interval, and we get $|f (x) - f
(a) | \leqslant M|x - a|$ from the rule of bounded change. As we will see
later (for general moduli of continuity), the analogy runs even deeper, and in
fact differentiation of ULD functions is related to factoring in the class of
Lipschitz functions the same way as differentiation of polynomials is related
to their factoring.

\subsubsection{A comparison with the non-standard \ analysis approach}

In non-standard analysis the derivatives of functions differentiable on a
hyperreal interval are automatically continuous, the proof goes the same way,
except we say that $\frac{f (x) - f (a)}{x - a} - f' (a)$ and $\frac{f (a) - f
(x)}{a - x} - f' (x)$ are infinitely small when $x - a$ is, and conclude that
in this case $f' (x) - f' (a)$ is infinitely small too. It is this fact that
makes the non-standard approach to calculus simple. More generally, many
pointwise estimates \ on a hyperreal interval are in fact uniform. In uniform
differentiation theory we work with uniform estimates directly and get the
results much cheaper, without any infinitesimals that are not constructive.
See http://en.wikipedia.org/wiki/Hyperreal\_number where I wrote a section
``An intuitive approach to the ultrapower construction'' and references there.

\subsection{Integration, existence of a primitive and Newton-Leibniz}

It is easy to integrate polynomials and rational functions since
antiderivatives can be written down explicitly in terms of the elementary
functions, but this situation is rather exceptional. Now, we know that the
derivative of any ULD function is Lipschiz, and we ask if an antiderivative
exists for any Lipschitz function, in what sense it exists, and how it can be
calculated. The idea is to define the definite integral as the area under the
graph and then to make sense out of the notion of the area by constructing
explicit approximations (pretty much following the approach of the Greeks,
later developed by Riemann, Darboux, Jordan, and Lebesque), and then prove
Newton-Leibniz. The case of Lipschitz, and other uniformly continuous
functions, is particularly simple, and requires hardly any sophistication. A
picture (that is worth a 1000 words) is available on page 13 at
http://www.mathfoolery.org/talk-2004.pdf and pages 43-44 at
http://www.mathfoolery.org/lathead.pdf with a proof of Newton-Leibniz.

\section{Other Moduli of Continuity}

Sometimes calculus based on Lipschitz estimates is too restrictive, for
example, the function $x^{3 / 2}$ has $\sqrt{x}$ for the derivative, which is
not Lipschiz, since it grows too fast near $x = 0$. To treat this function as
differentiable, we can relax the estimate (2) defining differentiability to \
$|f (x) - f (a) - f' (a) (x - a) | \leqslant K|x - a|^{3 / 2}$. More
generally, we can use the inequality
\[ |f (x) - f (a) - f' (a) (x - a) | \leqslant K|x - a|m (|x - a|) \]
with some \tmtextsl{modulus of continuity} $m$ to define
$m$\tmtextsl{-differentiability}, $m (x) = \sqrt{x}$ is an example, for any
positive $\gamma \leqslant 1$, $x^{\gamma}$ is a more general example, the
corresponding differentiability is called \tmtextsl{uniform Holder}, with the
exponent $\gamma$ and the corresponding derivatives are \tmtextsl{Holder},
i.e., $|f' (x) - f' (a) | \leqslant H|x - a|^{\gamma}$ holds. In general, we
want $m$ to be defined for $x \geqslant 0,$ an increasing, continuous at $0$,
$m (0) = 0,$ and \tmtextsl{subadditive}, i.e., $m (x + y) \leqslant m (x) + m
(y)$. All the Lipschitz theory extends to the general moduli of continuity
with some obvious modifications, the derivatives are
$m$\tmtextsl{}-continuous, i.e., $|f' (x) - f' (a) | \leqslant K m (|x - a|)$
etc.

\subsection{An estimate of the difference quotient}

Let $m$ be a subadditive modulus of continuity, in particular, $m$ is
increasing, defined for $x \geqslant 0$, and $m (x) / x$ is decreasing for $x
> 0$, and let $f$ be a uniformly $m$ - differentiable function, i.e. there is
a uniform in $x$ and $a$ estimate with some constant $K$:

\begin{equation}
  |f (x) - f (a) - f' (a) (x - a) | \leqslant K|x - a|m (|x - a|)
\end{equation}
Let the difference quotient for $f$ be the 2-variable function
\[ Q_f (x, a) = (f (x) - f (a)) / (x - a) \tmop{for} x \neq a \tmop{and} Q_f
   (x, x) = f' (x) . \]
We want to demonstrate the inequality
\begin{equation}
  |Q_f (x, a) = Q_f (y, a) | \leqslant 2 K m (|x - y|),
\end{equation}
\[  \]
that means that the difference quotient is a uniformly $m$ - continuous. That
will justify the idea that uniform differentiation is factoring in the class
of $m$ - continuous functions of 2 variables.

Because only the increments of the independent variable and the corresponding
increments of the values of $f$ are involved in the difference quotient, we
can assume $a = 0 = f (0)$ and the inequality we want becomes
\begin{equation}
  |f (x) / x - f (y) / y| \leqslant 2 K m (|x - y|) .
\end{equation}
The case $x < 0 < y$ is easy because $|f (x) / x - f' (0) | \leqslant K m
(|x|) \tmop{and} |f (y) / y - f' (0) | \leqslant K m (|y|)$, so $|f (x) / x -
f (y) / y) | \leqslant K (m (|x|) + m (|y|)) \leqslant 2 K m (|x - y|)$
because $m$ is increasing.

The case of $x$ and $y$ of the same sign, say, $0 < x < y$ is a bit more
delicate. First we notice that adding any linear function to $f$ does not
change $f (x) / x - f (y) / y$, so we can assume that $f' (x) = 0$. The
left-hand side of the inequality we want to prove can be rewritten as $| ((y -
x) f (x) - x (f (y) - f (x))) / (x y) |$. Now, $|f (y) - f (x) | \leqslant K
(y - x) m (y - x)$ because $f' (x) = 0$, and also $|f (x) | \leqslant K x m
(x)$ because $f (0) = 0$. So it is enough to show that $\frac{y - x}{y} m (x)
\leqslant m (y - x)$. Again the case $x \leqslant y - x$ is easy because $m$
is increasing. We only have to use subadditivity of $m$ when $y - x \leqslant
x$. In this case $m (x) / y \leqslant m (x) / x \leqslant m (y - x) / (y - x)$
and we are done.

\subsection{Epsilon-delta and moduli of continuity}

We used different moduli of continuity to describe uniform continuity and
differentiability. The question is: ``how much of the classical theory of
continuous and smooth functions do we miss, if any?'' The answer to this
question is ``nothing.'' Let us consider uniform continuity, uniform
differentiability is analogous.

The classical way to describe uniform continuity of a function $f$ is to say
that for any $\varepsilon > 0$ there is $\delta > 0$ such that $|f (x) - f (a)
| < \varepsilon$ when $|x - a| < \delta$.

We want to show that there is a modulus of continuity $m$, such that the
inequality $|f (x) - f (a) | \leqslant m (|x - a|)$ holds. Let us consider the
following function: $g (h) = \sup \{|f (x) - f (a) | : |x - a| \leqslant h\}.$
We know that $g$ will be positive, increasing, and $g (h) \longrightarrow 0$
as $h \longrightarrow 0$, so $g$ will become continuous at $0$ if we put $g
(0) = 0$. Now, on $(h, y)$ plane consider the set $\{(y, h) : y \leqslant g
(h)\}$ of points under the graph of $g$. Take the convex hull of this set. The
upper edge of this convex hull will be the graph of a concave (and therefore
subadditive) modulus of continuity for $f$.

It is needless to say that in some questions (such as topological
classification of dynamical systems) keeping track of the particular moduli of
continuity may be a nuisance, and not fruitful. Then we can throw all the
uniformly continuous or uniformly differentiable functions into one big pile
and enjoy the generality. \

\section{Some Pedagogical implications}

\subsection{Calculus by problem solving: a still unrealized dream}

http://www.mathfoolery.org/Problem\_sets/hw.html

\section{Many Variables}

\subsection{Differentiability}

\ Similar to the case of one variable, we define differentiability by the
inequality
\[ |f (x + h) - f (x) - f' (x) h| \leqslant K|h|m (|h|) . \]
Here $|.|$ denotes some norm, for example, the Euclidean norm, $f' (x)$ is a
linear map depending on $x$, $K$ is a constant and $m$ is a modulus of
continuity.

\subsection{Automatic continuity of the derivative}

We want to show that the uniform derivative is uniformly continuous with the
modulus of continuity $m$ from the definition, i.e., the inequality
\[ |f' (x + h) - f' (x) | \leqslant L m (|h|) \]
holds for some constant $L$ that will depend on $K$ in the definition. Here
$|.|$ is the norm of the linear operators, $|A | = \sup \{|A k|, |k| = 1\}.$

The idea of the simplest proof I could come up with is the following. There
are two ways to get from $x$ to $x + h + k.$ We can go directly, or we can go
from $x$ to $x + h$ first and then from $x + h$ to $x + h + k$. The
corresponding increments of the function $f$ should be the same. Now consider
the approximation of these increments by the differentials.
\[ |f (x + h) - f (x) - f' (x) h| \leqslant K|h| m (|h|) \]
\[ |f (x + h + k) - f (x + h) - f' (x + h) k| \leqslant K|k| m (|k|) \]
\[ | - f (x + h + k) + f (x) + f' (x) (h + k) | \leqslant K|m + k| m (|h + k|)
\]
By ``adding'' all of these inequalities and using the triangle inequality, $|a
+ b| \leqslant |a| + |b|$, and linearity, $f' (x) (h + k) = f' (x) h + f' (x)
k,$ we conclude that
\[ |f' (x) k - f' (x + h) k| \leqslant K (|h|m (|h|) + |k|m (|k|) + |h + k|m
   (|h + k|)) . \]
But $|h + k| \leqslant |h| + |k|$ and $m (|h + k|) \leqslant m (|h| + |k|)
\leqslant m (|h|) + m (|k|)$ (triangle, $m$ is increasing and subadditive).
Finally, by taking $|k| = |h|$, we get
\[ | (f' (x + h) - f' (x)) k| = |f' (x) k - f' (x + h) k| \leqslant 6 K m
   (|h|) |k| \]
that means that $|f' (x + h) - f' (x) | \leqslant^{} 6 K m (|h|)$, so we can
take $L = 6 K.$ Done.

\subsection{The equality of the mixed derivatives}

Probably the simplest way to understand why $f_{x y} = f_{y x}$ is to use the
Green's formula. Here is how. Let us consider a rectangle $A B C D$ on the
$(x, y) -$plane where $f$ is defined. $A = (a, c),$ $B = (b, c),$ $C = (b, d)$
and $D = (a, d) .$ There are two ways to get from $A$ to $C$. We can go from
$A$ to $B$ and then from $B$ to $C$, or we can go from $A$ to $D$ and then
from $D$ to $C$. The total change in $f$ should be the same for both ways. Let
us write down this change in terms of the line integrals of the partial
derivatives.$\text{}$
\[ \text{$f (A) - f (C) = f (B) - f (A) + f (C) - f (B) = \int_a^b f_x (x, c)
   d x + \int_c^d f_y (b, y) d y \text{}$} \]
and on the other hand,
\[ f (A) - f (C) = f (D) - f (A) + f (C) - f (D) = \int_{c^{}}^d f_y (a, y) d
   y + \int_a^b f_x (x, d) d x \]
so we have
\[ \int_a^b (f_x (x, d) - f_x (x, c)) d x -_{}^{} \int_c^d (f_y (b, y) - f_y
   (a, y)) d y = 0, \]
but $f_x (x, d) - f_x (x, c) = \int_c^d f_{x y} (x, y) d y$ and $f_y (b, y) -
f_y (a, y) = \int_a^b f_{y x} (x, y) d x$ and, replacing the iterated
integrals with the double integrals, we conclude that$\int_{A B C D} (f_{x y}
- f_{y x}) d x d y$ = 0 for any rectangle $A B C D$. It is only possible if
$f_{x y} - f_{y x} = 0$, so $f_{x y} = f_{y x} $ and we are done.

\end{document}